# Three-dimensional Alpha Shapes[1]

Herbert Edelsbrunner[2] and Ernst P Mücke[2]

**Abstract.** *Frequently, data in scientific computing is in its abstract form a finite point set in space, and it is sometimes useful or required to compute what one might call the "shape" of the set. For that purpose, this paper introduces the formal notion of the family of $\alpha$-shapes of a finite point set in $I\!R^3$. Each shape is a well-defined polytope, derived from the Delaunay triangulation of the point set, with a parameter $\alpha \in I\!R$ controlling the desired level of detail. An algorithm is presented that constructs the entire family of shapes for a given set of size n in time $O(n^2)$, worst case. A robust implementation of the algorithm is discussed and several applications in the area of scientific computing are mentioned.*

**Categories and Subject Descriptors.** F.2.2 [**Analysis of Algorithms and Problem Complexity**]: Nonnumerical algorithms and problems — *geometrical problems and computations*; G.4 [**Mathematics of Computing**]: Mathematical Software — *reliability and robustness*; I.2.10 [**Artificial Intelligence**]: Vision and Scene Understanding — *representations, data structures, and transforms*; *shape*; I.3.4 [**Computer Graphics**]: Graphics utilities — *application packages*; *graphics packages*; I.3.5 [**Computer Graphics**]: Computational Geometry and Object Modeling — *curve, surface, solid, and object representations*; *geometric algorithms, languages, and systems*; J.2 [**Computer Applications**]: Physical sciences and engineering.

ACM Transactions on Graphics, 13(1):43–72, 1994.

---

[1]Research of both authors is supported by the National Science Foundation under grant CCR-8921421. Work by the first author is also supported by the Alan T Waterman award, grant CCR-9118874. Any opinions, findings, conclusions, or recommendations expressed in this publication are those of the authors and do not necessarily reflect the view of the National Science Foundation.

[2]Address: Department of Computer Science, University of Illinois at Urbana-Champaign, 1304 W Springfield Ave, Urbana, IL 61801, USA.

# 1  Introduction

The geometric notion of "shape" has no associated formal meaning. This is in sharp contrast to other geometric notions, such as diameter, volume, convex hull, etc. The goal of this paper is to offer a concrete and formal definition of shape that can be computed and applied. It is not supposed to possibly cover the entire range of meanings the term "shape" carries in our contemporary language, even if restricted to geometric contexts. Nevertheless, it is sufficiently flexible to facilitate a wide range of applications including the study of molecular structures and the distribution of galaxies in our universe (see section 7).

More specifically, the topic of this paper is the definition and computation of the shape of a finite point set in three-dimensional Euclidean space, $I\!\!R^3$. Intuitively, we think of the set as a cloud of points and we talk about the shape of this cloud. A peculiar aspect of the common usage of the word "shape" is that its meaning varies with the degree of detail intended. This aspect will be reflected by defining a one-parametric family of shapes ranging from "fine" and "local" to "crude" and "global".

A fair amount of related work has been done for point sets in $I\!\!R^2$, and some for point sets in $I\!\!R^3$. Jarvis [23] was one of the first to consider the problem of computing the shape as a generalization of the convex hull of a planar point set. A mathematically rigorous definition of shape was later introduced by Edelsbrunner, Kirkpatrick, and Seidel [14]. Their notion of $\alpha$-shapes is the two-dimensional analogue of the spatial notion described in this paper. Two-dimensional $\alpha$-shapes are related to the dot patterns of Fairfield [17, 18] and the circle diagrams used in bivariate cluster analysis (see for example Moss [33]). Different graph structures that serve similar purposes are the Gabriel graph [31], the relative neighborhood graph [39], and their parameterized version, the $\beta$-skeleton [26].

For $I\!\!R^3$, Boissonnat [1] suggested the use of Delaunay triangulations in connection with heuristics to "sculpture" a single connected shape of a point set. Our concept of shape is more general and mathematically well defined. More recently, Veltkamp [40] also generalized the above-mentioned two-dimensional graph structures to three-dimensions, calling them $\gamma$-graphs. Finally, note the superficial similarity between $\alpha$-shapes and isosurfaces in $I\!\!R^3$. The latter is a popular concept in volume visualization (see for example [8, 29]).

The outline of this paper is as follows. A formal definition of $\alpha$-shapes, along with illustrations, is presented in section 2. Geometric concepts related to $\alpha$-shapes are discussed in section 3. These are $\alpha$-hulls, $\alpha$-diagrams (also known as space-filling diagrams), $\alpha$-complexes, Delaunay triangulations, and Voronoi diagrams. A combinatorial analysis of $\alpha$-shapes is offered in section 4. Using Delaunay triangulations, it is fairly easy to compute $\alpha$-shapes in $I\!\!R^3$. The resulting algorithm is sketched in section 5. Given a set of $n$ points in $I\!\!R^3$, it constructs a convenient implicit representation of the family of *all* $\alpha$-shapes in time $O(n^2)$, worst case. This algorithm has been implemented by the second author of this paper. In section 6 we report on essential aspects of the implementation, such as its data structures, how it achieves robustness, and how it performs in practice. Section 7 discusses some application problems that might benefit from the use of $\alpha$-shapes. Finally, section 8 considers possible extensions of the material presented in this paper.

# 2  Alpha Shapes in Space

This section gives an intuitive description as well as a formal definition of three-dimensional $\alpha$-shapes. Both are supported by illustrations that show point sets with sample members of their $\alpha$-shape family. The beauty and elegance of the concept of an $\alpha$-shape will hopefully be obvious after reading section 3 where relationships to other natural geometric concepts are revealed.

**Intuitive Description.** Conceptually, $\alpha$-shapes are a generalization of the convex hull of a point set. Let $S$ be a finite set in $\mathbb{R}^3$ and $\alpha$ a real number with $0 \leq \alpha \leq \infty$. The $\alpha$-shape of $S$ is a polytope that is neither necessarily convex nor necessarily connected. For $\alpha = \infty$, the $\alpha$-shape is identical to the convex hull of $S$. However, as $\alpha$ decreases, the $\alpha$-shape shrinks by gradually developing cavities. These cavities may join to form tunnels, and even holes may appear (see figure 1).

Intuitively, a piece of the polytope disappears when $\alpha$ becomes small enough so that a sphere with radius $\alpha$, or several such spheres, can occupy its space without enclosing any of the points of $S$. Think of $\mathbb{R}^3$ filled with styrofoam and the points of $S$ made of more solid material, such as rock. Now imagine a spherical eraser with radius $\alpha$. It is omnipresent in the sense that it carves out styrofoam at all positions where it does not enclose any of the sprinkled rocks, that is, points of $S$. The resulting object will be called the $\alpha$-hull (see section 3). To make things more feasible we straighten the surface of the object by substituting straight edges for the circular ones and triangles for the spherical caps. The obtained object is the $\alpha$-shape of $S$ (see figure 2). It is a polytope in a fairly general sense: it can be concave and even disconnected, it can contain two-dimensional patches of triangles and one-dimensional strings of edges, and its components can be as small as single points. The parameter $\alpha$ controls the maximum "curvature" of any cavity of the polytope.

**General Position.** Throughout this paper we assume that the points of $S$ are in general position. For the time being, this means that no 4 points lie on a common plane; no 5 points lie on a common sphere; and for any fixed $\alpha$, the smallest sphere through any 2, 3, or 4 points of $S$ has a radius different from $\alpha$. The general-position assumption will later be extended when convenient (see section 5.3). It simplifies forthcoming definitions, discussions, and algorithms, and is justified by a programming technique known as SoS [15]. This method simulates an infinitesimal perturbation of the points on the level of geometric predicates and relieves the programmer from the otherwise necessary case analysis (see section 6.2).

**Formal Definition.** For $0 < \alpha < \infty$, let an $\alpha$-*ball* be an open ball with radius $\alpha$. For completeness, a 0-*ball* is a point and an $\infty$-*ball* is an open half-space. An $\alpha$-ball $b$ is *empty* if $b \cap S = \emptyset$. Any subset $T \subseteq S$ of size $|T| = k+1$, with $0 \leq k \leq 3$, defines a $k$-*simplex* $\sigma_T$ that is the convex hull of $T$, also denoted by $\mathrm{conv}(T)$. The general-position assumption assures that all $k$-simplices are properly $k$-dimensional. For $0 \leq k \leq 2$, a $k$-simplex $\sigma_T$ is said to be $\alpha$-*exposed* if there is an empty $\alpha$-ball $b$ with $T = \partial b \cap S$, where $\partial b$ is the sphere or plane bounding $b$. A fixed $\alpha$ thus defines sets $F_{k,\alpha}$ of $\alpha$-exposed $k$-simplices for $0 \leq k \leq 2$. The $\alpha$-*shape* of $S$, denoted by $\mathcal{S}_\alpha$, is the polytope whose boundary consists of the triangles in $F_{2,\alpha}$, the edges in $F_{1,\alpha}$, and the vertices in $F_{0,\alpha}$ (see figures 1 and 2). The $k$-simplices in $F_{k,\alpha}$ are also called the $k$-*faces* of $\mathcal{S}_\alpha$.

We still need to specify which connected components of $\mathbb{R}^3 - \partial \mathcal{S}_\alpha$ are interior and which are exterior to $\mathcal{S}_\alpha$. Fix the value of $\alpha$ and notice that for each $\alpha$-exposed triangle $\sigma_T$ there are two (not necessarily empty) $\alpha$-balls, $b_1 \neq b_2$, so that $T \subseteq \partial b_1$ and $T \subseteq \partial b_2$. If both $\alpha$-balls are empty then $\sigma_T$ does not belong to the boundary of the interior of $\mathcal{S}_\alpha$. Otherwise, assume that $b_1$ is empty and that $b_2$ is not. In this case, $\sigma_T$ bounds the interior of $\mathcal{S}_\alpha$. More specifically, the interior of $\mathcal{S}_\alpha$ and the center of $b_1$ lie on different sides of $\sigma_T$. The definition of interior and exterior of $\mathcal{S}_\alpha$ is possibly more natural in the context of Delaunay triangulations and $\alpha$-complexes as described in section 3.

Figure 1: Two tori. [$n = 800$, $|F_2| = 12197$]

The points are randomly generated on the surface of two linked tori. Six different $\alpha$-shapes for values of $\alpha$ decreasing from top to bottom and left to right are shown. The first shape is the convex hull, for $\alpha = +\infty$; the last shape is the point set itself, for $\alpha = 0$. The $\alpha$-value used in the forth frame neatly separates the two tori. Further decreasing $\alpha$ disassembles the shape. Singular triangles, which do not bound the interior of the shape, are shown in darker color.

Figure 2: Bust. [$n = 2630$, $|F_2| = 35196$]

This point set is based on a demo data set for Silicon Graphics' Solidview program, of course, without any of the original connectivity information. The erasing sphere is shown to the right of the shape. Apart from a dense conglomerate of points representing part of the person's brain and brain stem, the set is basically hollow with most points representing skin.

# 3  Related Geometric Concepts

There are quite a few natural geometric concepts that are closely related to $\alpha$-shapes. Some of them are discussed in this section. In each case, the emphasis is on how the concept is related to $\alpha$-shapes and how this relation can enrich our understanding of $\alpha$-shapes. Section 3.1 discusses $\alpha$-hulls and $\alpha$-diagrams. Section 3.2 briefly reviews Delaunay triangulations and their dual incarnations known as Voronoi diagrams. The relevance of the Delaunay triangulation of a point set is that each $\alpha$-shape of the set is the underlying space of a subcomplex of the triangulation. These subcomplexes are termed $\alpha$-complexes and defined in section 3.3. Extensions of these notions are mentioned in section 3.4.

## 3.1  Alpha Hulls and Alpha Diagrams

Recall from section 2 that for positive real $\alpha$ an $\alpha$-ball is defined as an open ball with radius $\alpha$. For $\alpha = 0$, it is a point, and for $\alpha = \infty$, it is an open half-space. Given a finite point set $S \subseteq \mathbb{R}^3$, an $\alpha$-ball is empty if $b \cap S = \emptyset$. For $0 \leq \alpha \leq \infty$, the $\alpha$-*hull* of $S$, denoted by $\mathcal{H}_\alpha$, is defined as the complement of the union of all empty $\alpha$-balls. This is the formal counterpart of the styrofoam object described in section 2. Sample members of the continuous family of $\alpha$-hulls are the convex hull of $S$, for $\alpha = \infty$, and $S$ itself, for $\alpha$ sufficiently small. Observe that $\mathcal{H}_{\alpha_1} \subseteq \mathcal{H}_{\alpha_2}$ if $\alpha_1 \leq \alpha_2$.

Another interesting concept defined by $\alpha$-balls is what we call the $\alpha$-*diagram* of $S$, denoted by $\mathcal{U}_\alpha$. For $0 < \alpha < \infty$, $\mathcal{U}_\alpha$ is the union of all $\alpha$-balls whose centers are points in $S$.[3] Observe that a point $x \in \mathbb{R}^3$ belongs to $\mathcal{U}_\alpha$ iff the $\alpha$-ball centered at $x$ is not empty. Denote this $\alpha$-ball by $b_x$. This implies the following close relationship between $\mathcal{H}_\alpha$ and $\mathcal{U}_\alpha$.

$$x \in \mathcal{U}_\alpha \iff b_x \cap \mathcal{H}_\alpha \neq \emptyset, \text{ and}$$
$$x \in \mathcal{H}_\alpha \iff b_x \subseteq \mathcal{U}_\alpha.$$

Consider the boundary of $\mathcal{U}_\alpha$. It consists of spherical caps, circular arcs, and vertices which we call *corners*. These are the 2-, 1-, and 0-faces of $\mathcal{U}_\alpha$. These caps, arcs, and corners are in close correspondence with the vertices, edges, and triangles of $\mathcal{S}_\alpha$. Some definitions are needed to describe this correspondence.

Let $\alpha$ be fixed, with $0 < \alpha < \infty$, and let $T$ be a subset of $S$ of size $|T| = k + 1$, with $0 \leq k \leq 2$. Define $K_T = \bigcap_{p \in T} \partial b_p$, where $b_p$ is the $\alpha$-ball centered at $p$, as before. Besides general position assume $K_T \neq \emptyset$. For $|T| = 1$, $K_T$ is a sphere; for $|T| = 2$, $K_T$ is a circle; and for $|T| = 3$, $K_T$ is a pair of points. It follows from the definitions that $T$ is $\alpha$-exposed iff $K_T$ contains at least one face of the boundary of $\mathcal{U}_\alpha$. If $|T| = 1$ this face is a cap; if $|T| = 2$ it is an arc; and if $|T| = 3$ it is a corner. This fact can be expressed as follows.

$$\sigma_T \text{ is a } k\text{-face of } \mathcal{S}_\alpha, \text{ with } 0 \leq k \leq 2 \iff K_T \text{ contains at least one } (2-k)\text{-face of } \mathcal{U}_\alpha.$$

---

[3] In chemistry and biology, $\alpha$-diagrams are known as space-filling diagrams. However, there they are usually not restricted to equally large balls. This restriction can be removed with weighted $\alpha$-shapes and $\alpha$-diagrams (see section 3.4).

Moreover, the number and structure of $(2-k)$-faces contained in $K_T$ are reflected by $\sigma_T$ and the way it is embedded in $\mathcal{S}_\alpha$. For example, if $|T| = 3$, then $K_T$ contains no, one, or two corners of $\mathcal{U}_\alpha$. First, it contains no corner iff $\sigma_T$ is not a triangle of $\mathcal{S}_\alpha$. This is mentioned above. Second, $K_T$ contains one corner iff $\sigma_T$ bounds the interior of $\mathcal{S}_\alpha$. Third, both points of $K_T$ are corners of $\mathcal{U}_\alpha$ iff both sides of $\sigma_T$ face the outside of $\mathcal{S}_\alpha$, that is, $\sigma_T$ is a triangle of $\mathcal{S}_\alpha$ that does not bound its interior. Such a triangle will be called singular in section 5.2. Similarly, if $|T| = 2$ and $\sigma_T$ is an edge of $\mathcal{S}_\alpha$, then every angular interval between two incident triangles of $\mathcal{S}_\alpha$ that faces the outside corresponds to an arc of $\mathcal{U}_\alpha$ contained in $K_T$. If there is no incident triangle (in this case, $\sigma_T$ is an singular edge) then the entire 1-sphere $K_T$ belongs to the boundary of $\mathcal{U}_\alpha$. Analogous statements can be made for $|T| = 1$ where spherical caps on $K_T$ correspond to solid angles around the vertex $\sigma_T$ that face the outside of $\mathcal{S}_\alpha$.

All this amounts to a one-to-one correspondence between the $(2-k)$-faces of $\mathcal{U}_\alpha$ and the $k$-faces of $\mathcal{S}_\alpha$, provided the latter are interpreted with multiplicities reflecting the number of exposed sides, angular intervals, or solid angles. The one-to-one correspondence also preserves incidences. This suggests that $\mathcal{U}_\alpha$, or more specifically, the boundary of $\mathcal{U}_\alpha$, be represented by $\mathcal{S}_\alpha$, or more specifically, the faces of $\mathcal{S}_\alpha$ and their incidences.

## 3.2 Delaunay Triangulations and Voronoi Diagrams

A finite point set $S \subseteq \mathbb{R}^3$ defines a special triangulation known as the Delaunay triangulation of $S$ (see for example [11, 35]). Assuming general position of the points, this triangulation is unique and decomposes the convex hull of $S$ into tetrahedra. The triangulation is named after the Russian geometer Boris Delaunay (also Delone) who introduced it in his seminal paper [6] in 1934. As explained below, the Delaunay triangulation of $S$ is dual to another complex defined by $S$ known as the Voronoi diagram [41, 42]. Both complexes are related in an interesting way to the family of all $\alpha$-shapes of $S$. The relationship between $\alpha$-shapes and Delaunay triangulations will be of particular importance for this paper.

**Delaunay Triangulations.** For $0 \leq k \leq 3$, let $F_k$ be the set of $k$-simplices $\sigma_T = \text{conv}(T)$, $T \subseteq S$ and $|T| = k+1$, for which there are empty open balls $b$ with $\partial b \cap S = T$. Notice that $F_0 = S$. The *Delaunay triangulation* of $S$, denoted by $\mathcal{D}$, is the simplicial complex defined by the tetrahedra in $F_3$, the triangles in $F_2$, the edges in $F_1$, and the vertices in $F_0$. By definition, for each simplex $\sigma_T \in \mathcal{D}$, there exist values of $\alpha \geq 0$ so that $\sigma_T$ is $\alpha$-exposed. Conversely, every face of $\mathcal{S}_\alpha$ is a simplex of $\mathcal{D}$. This implies the following relationship between the Delaunay triangulation and the boundary of $\mathcal{S}_\alpha$.

$$F_k = \bigcup_{0 \leq \alpha \leq \infty} F_{k,\alpha}, \text{ for } 0 \leq k \leq 2.$$

We take advantage of this relationship by representing the family of $\alpha$-shapes of $S$ implicitly by the Delaunay triangulation of $S$. This will be described in detail in section 5.

**Voronoi Diagrams.** For a point $p \in S$, define $V(p)$, the *Voronoi cell* of $p$, as the set of points $x \in \mathbb{R}^3$ so that the Euclidean distance between $x$ and $p$ is less than or equal to the distance between $x$ and any other point of $S$. Each Voronoi cell is a convex polyhedron, and the collection of all Voronoi cells, one for each point of $S$, defines the *Voronoi diagram* of $S$, denoted by $\mathcal{V}$. We call a Voronoi cell also a 3-*cell* of $\mathcal{V}$. Each 2-*cell* of $\mathcal{V}$ is the intersection of two Voronoi cells; each 1-*cell* or edge is the intersection of three 3-cells; and each 0-*cell* or vertex is the intersection of four 3-cells. There is a natural one-to-one correspondence between the $k$-simplices of $\mathcal{D}$ and the $(3-k)$-cells of $\mathcal{V}$. Let $T$ be a subset of $S$ of size $|T| = k+1$, with $0 \leq k \leq 3$, and define $V_T = \bigcap_{p \in T} V(p)$.

$$\sigma_T \text{ is a } k\text{-simplex of } \mathcal{D} \iff V_T \text{ is a } (3-k)\text{-cell of } \mathcal{V}, \text{ for } 0 \leq k \leq 3.$$

This correspondence preserves (or reverses) incidences which implies that $\mathcal{D}$ and $\mathcal{V}$ are indeed dual to each other.

Observe that $V_T$ is the set of all points $x$ for which there exists an empty open ball $b_x$ centered at $x$ with $T \subseteq \partial b_x \cap S$. Equality holds iff $x$ belongs to the relative interior of $V_T$. It follows that $\sigma_T$ is $\alpha$-exposed iff there is a point $x$ in the relative interior of $V_T$ whose distance from the points in $T$ is $\alpha$. Since $V_T$ is convex there is a single interval so that $\sigma_T$ is $\alpha$-exposed iff $\alpha$ belongs to this interval. We will exploit this fact later when we discuss how to decide when a simplex of $\mathcal{D}$ is $\alpha$-exposed.

## 3.3 Alpha Complexes

Since all faces of $\mathcal{S}_\alpha$ are simplices of $\mathcal{D}$, it follows that the interior of $\mathcal{S}_\alpha$ is naturally triangulated by the tetrahedra of $\mathcal{D}$. This idea leads to the concept of $\alpha$-complexes as defined shortly. A (three-dimensional) *simplicial complex* is a collection $\mathcal{C}$ of closed $k$-simplices, for $0 \leq k \leq 3$, that satisfies the following two properties.

(i) If $\sigma_T \in \mathcal{C}$ then $\sigma_{T'} \in \mathcal{C}$ for every $T' \subseteq T$. In other words, with every simplex $\sigma_T$, $\mathcal{C}$ contains all faces of $\sigma_T$ as well.

(ii) If $\sigma_T, \sigma_{T'} \in \mathcal{C}$, then either $\sigma_T \cap \sigma_{T'} = \emptyset$, or $\sigma_T \cap \sigma_{T'} = \sigma_{T \cap T'} = \operatorname{conv}(T \cap T')$. Note that (i) implies that this face is also in $\mathcal{C}$. In other words, the intersection of any two simplices in $\mathcal{C}$ is either empty or a face of both.

A subset $\mathcal{C}' \subseteq \mathcal{C}$ is a *subcomplex* of $\mathcal{C}$ if it is also a simplicial complex.

Each $k$-simplex $\sigma_T$ of $\mathcal{D}$ defines an open ball $b_T$ bounded by the smallest sphere $\partial b_T$ that contains all points of $T$. Let $\varrho_T$ be the radius of $b_T$. For $k = 3$, $\partial b_T$ is the circumsphere of $\sigma_T$; for $k = 2$, the circumcircle of $\sigma_T$ is a great circle of $\partial b_T$; and for $k = 1$, the two points in $T$ are antipodal on $\partial b_T$. Call $\partial b_T$ the *smallest circumsphere* and $\varrho_T$ the *radius* of $\sigma_T$. For $1 \leq k \leq 3$ and $0 \leq \alpha \leq \infty$, define $G_{k,\alpha}$ as the set of $k$-simplices $\sigma_T \in \mathcal{D}$ for which $b_T$ is empty and $\varrho_T < \alpha$. Furthermore, define $G_{0,\alpha} = S$, for all $\alpha$. The sets $G_{k,\alpha}$ do not necessarily define a simplicial complex because it can happen that $G_{3,\alpha}$ contains a tetrahedron but not all triangles of this tetrahedron belong to $G_{2,\alpha}$. Similarly for triangles and edges. With this in mind, we define the $\alpha$-*complex* of $S$, denoted by $\mathcal{C}_\alpha$, as the simplicial complex whose $k$-simplices are either in $G_{k,\alpha}$ or they bound $(k+1)$-simplices of $\mathcal{C}_\alpha$. By definition, $\mathcal{C}_{\alpha_1}$ is a subcomplex of $\mathcal{C}_{\alpha_2}$ if $\alpha_1 \leq \alpha_2$.

The *underlying space* of $\mathcal{C}_\alpha$, denoted by $|\mathcal{C}_\alpha|$, is the union of all simplices of $\mathcal{C}_\alpha$, or in other words, the part of $I\!\!R^3$ covered by $\mathcal{C}_\alpha$. Thus, the underlying space of $\mathcal{C}_\alpha$ is a polytope in the sense specified in section 2. Indeed, we have the following most important property of $\mathcal{C}_\alpha$, which we present without proof.

For all $0 \leq \alpha \leq \infty$, $\mathcal{S}_\alpha = |\mathcal{C}_\alpha|$.

This can be considered an alternative definition of $\alpha$-shapes. It makes it easy to specify the intervals of the simplices of $\mathcal{D}$ alluded to in the above discussion of Voronoi diagrams. For example, let $\sigma_T$ be a $k$-simplex of $\mathcal{D}$. If $b_T$ is empty then $\sigma_T$ belongs to $\mathcal{C}_\alpha$ iff $\alpha \in (\varrho_T, \infty]$.[4]

## 3.4 Extensions

The definitions presented in sections 2 and 3 above can be extended in various ways. So far we refrained from mentioning these extensions in order to avoid unnecessary complications and to be faithful to the currently available implementation of the concepts in this paper. For completeness, negative values of $\alpha$, weighted points, and higher dimensions are now briefly discussed.

**Negative Alpha Values.** This extension has been described in [14] for the two-dimensional case.[5] For negative $\alpha$, $\alpha$-complexes are most naturally defined as subcomplexes of the so-called furthest-point Delaunay triangulation of $S$ (see for example [11, 35]). For $T \subseteq S$ and $|T| = 4$, the tetrahedron $\sigma_T$ belongs to this triangulation iff $b_T$ contains *all* points of $S - T$. The $\alpha$-shape is, again, the underlying space of the $\alpha$-complex. These shapes exhibit far less interesting geometric and topological properties than the ones for positive $\alpha$, and are thus less interesting for applications. We omit further details.

**Weighted Points.** Recall the relationship between $\alpha$-shapes and $\alpha$-diagrams described in section 3.1. It is interesting to consider diagrams for different ball sizes, and this is indeed done in chemistry and biology where space-filling diagrams are usually defined as unions of balls with arbitrary and thus possibly different radii. In order to represent such diagrams by polytopes similar to $\alpha$-shapes it is necessary to introduce weighted $\alpha$-shapes. These can be defined using subcomplexes of so-called regular triangulations (see for example [12, 28]). Given a finite set of points, each with a real weight, the regular triangulation is a unique simplicial complex whose underlying space is the convex hull of the point set. If all weights are the same then it equals the Delaunay triangulation of the points. Details can be found in [13].

**Higher Dimensions.** It is fairly straightforward to generalize all important concepts of this section (like $\alpha$-shapes, $\alpha$-hulls, $\alpha$-diagrams, $\alpha$-complexes, Delaunay triangulations, and Voronoi diagrams) to finite point sets $S$ in $\mathbb{R}^d$, for arbitrary dimension $d$. This generalization, combined with an extension to weighted points is developed in [13]. Note however, that the implementation details become progressively "hairier" with increasing dimension, and the worst-case complexity of the problem grows exponentially. For example, if $S$ is a set of $n$ points in $\mathbb{R}^d$ then the Delaunay triangulation of $S$ can consist of up to $\Theta(n^{\lfloor \frac{d+1}{2} \rfloor})$ faces (see [38]). Although the running time of the programs constructing $\alpha$-shapes will get substantially worse as $d$ increases, there might be applications that warrant implementations in low dimensions higher than three.

---

[4] Because of the general-position assumption we have $\alpha \neq \varrho_T$. It is therefore irrelevant whether the interval is open or closed at its left endpoint.

[5] The definitions in [14] are slightly different from the ones in this paper. In particular, their $\alpha$ is the same as minus one over $\alpha$ in this paper. Thus, our negative values of $\alpha$ correspond to their positive values.

# 4  Combinatorial Analysis

In contrast to $\alpha$-hulls and $\alpha$-diagrams, the $\alpha$-shapes of a finite point set form a discrete family, even though they are defined for *all* real numbers $\alpha$, with $0 \leq \alpha \leq \infty$. Indeed, $\mathcal{S}_{\alpha_1} \neq \mathcal{S}_{\alpha_2}$ iff $\bigcup_{k=0}^{3} G_{k,\alpha_1} \neq \bigcup_{k=0}^{3} G_{k,\alpha_2}$. Thus, $\mathcal{S}_{\alpha_1} \neq \mathcal{S}_{\alpha_2}$ iff there is an empty open ball bounded by a smallest circumsphere of an edge, triangle, or tetrahedron of $\mathcal{D}$ whose radius lies between $\alpha_1$ and $\alpha_2$. Such a radius is referred to as an $\alpha$-*threshold* because it separates two $\alpha$-shapes. The number of $\alpha$-shapes exceeds the number of $\alpha$-thresholds by one. It follows that one plus the total number of $k$-simplices of $\mathcal{D}$, for $1 \leq k \leq 3$, is an upper bound on the number of different $\alpha$-shapes. An upper bound on the number of simplices of $\mathcal{D}$ can be obtained using a relationship between Delaunay triangulations in $\mathbb{R}^3$ and certain convex polytopes in $\mathbb{R}^4$. We briefly describe this relationship in the next paragraph and refer to [11, 38] for details.

**Lifting Map.** Identify $\mathbb{R}^3$ with the $x_1 x_2 x_3$-space in $\mathbb{R}^4$, that is, the subspace $x_4 = 0$. The *lifting map* is a geometric transform that projects points $p = (\pi_1, \pi_2, \pi_3)$ in $\mathbb{R}^3$ along the $x_4$-axis onto the paraboloid of revolution $U: x_4 = \sum_{i=1}^{3} x_i^2$ in $\mathbb{R}^4$. Let $p_U = (\pi_1, \pi_2, \pi_3, \sum_{i=1}^{3} \pi_i^2)$ be the image of $p$, and define $S_U = \{p_U \mid p \in S\}$. The convex hull of $S_U$, $\mathrm{conv}(S_U)$, is a convex polytope with $n$ vertices in $\mathbb{R}^4$. A facet belongs to the *lower* boundary of this polytope if the polytope lies on the side of the positive $x_4$-axis of the hyperplane that contains the facet. Otherwise, it belongs to the *upper* boundary of $\mathrm{conv}(S_U)$. If we project all facets of the lower boundary of $\mathrm{conv}(S_U)$ parallel to the $x_4$-axis into $\mathbb{R}^3$, along with their subfaces, then we obtain the Delaunay triangulation of $S$ (see [11]).

**Upper Bounds.** According to the upper bound theorem for convex polytopes, the maximum numbers of 1-, 2-, and 3-faces of a convex polytope with $n \geq 5$ vertices in $\mathbb{R}^4$ are $\frac{1}{2}(n^2 - n)$, $n^2 - 3n$, and $\frac{1}{2}(n^2 - 3n)$, respectively (see for example [3, 32]). The lifting map implies the same upper bounds for $|F_k|$, with $1 \leq k \leq 3$. As a matter of fact, the upper bound for $|F_3|$ is one less than for the number of 3-faces of $\mathrm{conv}(S_U)$ because at least one 3-face belongs to the upper boundary of $\mathrm{conv}(S_U)$. By a result of [38], these bounds are tight even though the vertices of $\mathrm{conv}(S_U)$ are constrained to lie on a second-degree surface in $\mathbb{R}^4$. We summarize these results for $n = |S|$.

$$|F_0| = n, \quad |F_1| \leq \frac{1}{2}(n^2 - n), \quad |F_2| \leq n^2 - 3n, \quad \text{and} \quad |F_3| \leq \frac{1}{2}(n^2 - 3n - 2).$$

By adding one to the sum of the bounds for $|F_1|$, $|F_2|$, and $|F_3|$, we obtain the following result.

$S$ has at most $2n^2 - 5n$ different $\alpha$-shapes.

This bound is too pessimistic for two reasons. First, although the upper bounds on the number of simplices of $\mathcal{D}$ are tight, there are many fewer simplices for most point sets. For example, if the $n$ points are uniformly distributed in the unit ball then the expected number of simplices is only $O(n)$ (see Dwyer [9]). Second, not all edges and triangles of $\mathcal{D}$ have a smallest circumsphere that bounds an empty open ball. However, since the circumsphere of every tetrahedron of $\mathcal{D}$ bounds an empty open ball by definition, the number of different $\alpha$-shapes is always at least a fraction of the number of simplices of $\mathcal{D}$.

# 5 Algorithms

As described in section 3, the family of $\alpha$-shapes of a finite point set $S$ can be represented by the Delaunay triangulation of $S$. In this representation, each simplex of $\mathcal{D}$ is associated with an interval that specifies for which values of $\alpha$ the simplex belongs to the $\alpha$-shape. Section 5.1 discusses the construction of $\mathcal{D}$, and section 5.2 explains how the intervals of the simplices are computed. For completeness, section 5.3 gives the formulas that can be used to implement the required primitive operations.

## 5.1 Three-dimensional Delaunay Triangulations

The construction of Delaunay triangulations is a popular topic in the area of geometric algorithms (see for example [11, 35]). Indeed, various different approaches have been studied and described in the literature. Some approaches are based on the lifting map mentioned in section 4, which transforms the problem into one of constructing the convex hull of a four-dimensional point set.

The algorithm adopted for our implementation of $\alpha$-shapes has been suggested by Joe [25]. It is based on the idea of local transformations or flips.[6] The algorithm can be viewed as a generalization of the edge-flip method for two-dimensional triangulations by Lawson [27]. Note, however, that the straightforward generalization of the two-dimensional algorithm to $I\!\!R^3$ fails to always compute the Delaunay triangulation (see [24]). Nevertheless, the correctness of the flip algorithm in $I\!\!R^3$ can be established if the points are added one by one. We sketch the structure of the resulting incremental algorithm below and discuss the notion of a flip later.

**Incremental-flip Algorithm.**
    Sort $S$ along some fixed direction,
    and relabel the points so that $(p_1, p_2, \ldots, p_n)$ is the sorted sequence.
    Initialize $\mathcal{D}$ to the triangulation whose only tetrahedron is $\sigma_T$, $T = \{p_1, p_2, p_3, p_4\}$.
    **for** $i := 5$ **to** $n$ **do**
        Add $p_i$ by connecting it to all vertices, edges, and triangles of $\mathcal{D}$ visible from $p_i$.
        Use flips to transform $\mathcal{D}$ to the Delaunay triangulation of $\{p_1, p_2, \ldots, p_i\}$.
    **end for.**

Consider two tetrahedra $\sigma_{T'}$ and $\sigma_{T''}$ that share a common triangle $\sigma_T$. For example, $T = \{p_i, p_j, p_k\}$, $T' = T \cup \{p_u\}$, and $T'' = T \cup \{p_v\}$. Triangle $\sigma_T$ is called *locally Delaunay* if it belongs to the Delaunay triangulation of $T' \cup T''$. This is the case iff the point $p_v$ lies outside the sphere $\partial b_{T'}$. Local Delaunayhood is a necessary condition for $\sigma_T$ to belong to the Delaunay triangulation, but it is not sufficient. Nevertheless, Delaunay [6] proved that if *all* triangles are locally Delaunay then the triangulation is a Delaunay triangulation.

The incremental-flip algorithm needs to restore Delaunayhood whenever a new point is added to the triangulation. For this, it identifies triangles that are not locally Delaunay and tries to flip them. Let $\sigma_T$ be such a triangle bounding the tetrahedra $\sigma_{T'}$ and $\sigma_{T''}$. Again, assume $T = \{p_i, p_j, p_k\}$, $T' = T \cup \{p_u\}$, and $T'' = T \cup \{p_v\}$. If $\sigma_{T'} \cup \sigma_{T''}$ is convex then the triangle $\sigma_T$ can be replaced by the edge connecting $p_u$ and $p_v$. Together with this edge, the triangles connecting it with the three vertices of $T$ are added. This operation is called a *triangle-to-edge flip*. Otherwise, $\sigma_{T'} \cup \sigma_{T''}$ is not convex. Assume there is a third tetrahedron $\sigma_{T'''}$ that is spanned by four of the five points of $T' \cup T''$, for example, $T''' = \{p_i, p_j, p_u, p_v\}$. In this case, the three triangles incident to the edge connecting $p_i$ and $p_j$ can be replaced by a single triangle with endpoints $p_u$, $p_v$, and $p_k$. We call this an *edge-to-triangle flip*. If $\sigma_{T'''}$ does not exist then there is no flip that can remove $\sigma_T$.

---

[6] As recently shown, the flip algorithm can be extended to compute regular triangulations in arbitrary dimensions [16]. Regular triangulations are useful for weighted $\alpha$-shapes (see section 3.4 or [13]).

In spite of possible triangles that are neither locally Delaunay nor can be flipped according to the above rules, the correctness of the incremental-flip algorithm can still be established [25]. In the worst case, it takes time and storage $O(n^2)$, where $n = |S|$. Experiments provide evidence that it performs significantly better for most point sets. However, the worst case of $\Theta(n^2)$ cannot be avoided because $\mathcal{D}$ can have up to a quadratic number of simplices (see section 4).

## 5.2 Intervals and Face Classification

For each simplex $\sigma_T \in \mathcal{D}$ there is a single interval so that $\sigma_T$ is a face of the $\alpha$-shape $\mathcal{S}_\alpha$ iff $\alpha$ is contained in this interval. This was mentioned in section 3. It will be convenient to study these intervals for the $\alpha$-complex $\mathcal{C}_\alpha$ rather than the $\alpha$-shape. Also, we break each interval into three (possibly empty) parts that correspond to values of $\alpha$ for which the simplex is an interior, regular, or singular simplex of $\mathcal{C}_\alpha$.

A simplex $\sigma_T \in \mathcal{C}_\alpha$ is said to be

$$\begin{cases} interior & \text{if } \sigma_T \notin \partial \mathcal{S}_\alpha, \\ regular & \text{if } \sigma_T \in \partial \mathcal{S}_\alpha \text{ and it bounds some higher-dimensional simplex in } \mathcal{C}_\alpha, \text{ and} \\ singular & \text{if } \sigma_T \in \partial \mathcal{S}_\alpha \text{ and it does not bound any higher-dimensional simplex in } \mathcal{C}_\alpha. \end{cases}$$

Notice that there are Delaunay edges and triangles that can never be singular because their smallest circumsphere encloses other points of $S$.[7] Therefore, we call a simplex $\sigma_T \in \mathcal{D}$

$$\begin{cases} attached & \text{if } |T| = 2, 3 \text{ and } b_T \cap S \neq \emptyset, \text{ and} \\ unattached & \text{otherwise.} \end{cases}$$

Recall that $\varrho_T$ is the radius of the smallest circumsphere of $\sigma_T$. In order to break up the interval for which $\sigma_T$ belongs to $\mathcal{C}_\alpha$, we introduce values $\underline{\mu}_T$ and $\overline{\mu}_T$ for which $\sigma_T$ changes from singular to regular and from regular to interior, respectively. Before that, however, let up$(\sigma_T)$ be the set of all simplices in $\mathcal{D}$ that contain a simplex $\sigma_T \in \mathcal{D}$, with $|T| \leq 3$, as a proper face, that is,

$$\text{up}(\sigma_T) = \{\sigma_{T'} \in \mathcal{D} \mid T \subset T'\}.$$

---

[7]It is convenient to extend the general-position assumption so that no smallest circumsphere of two or three points of $S$ contains another point of $S$. A slightly more general assumption is the following. If a sphere contains three points of $S$ then no two of them are antipodal, and if it contains four points then no three lie on a great-circle of the sphere.

| | | singular | regular | interior |
|---|---|---|---|---|
| tetrahedron | | | | $(\varrho_T, \infty]$ |
| edge or triangle, | $\notin \partial\mathrm{conv}(S)$, unattached | $(\varrho_T, \underline{\mu}_T)$ | $(\underline{\mu}_T, \overline{\mu}_T)$ | $(\overline{\mu}_T, \infty]$ |
| | $\notin \partial\mathrm{conv}(S)$, attached | | $(\underline{\mu}_T, \overline{\mu}_T)$ | $(\overline{\mu}_T, \infty]$ |
| | $\in \partial\mathrm{conv}(S)$, unattached | $(\varrho_T, \underline{\mu}_T)$ | $(\underline{\mu}_T, \infty]$ | |
| | $\in \partial\mathrm{conv}(S)$, attached | | $(\underline{\mu}_T, \infty]$ | |
| vertex, | $\notin \partial\mathrm{conv}(S)$ | $[0, \underline{\mu}_T)$ | $(\underline{\mu}_T, \overline{\mu}_T)$ | $(\overline{\mu}_T, \infty]$ |
| | $\in \partial\mathrm{conv}(S)$ | $[0, \underline{\mu}_T)$ | $(\underline{\mu}_T, \infty]$ | |

Table 1: Intervals of $\alpha$ values for which $\sigma_T \in \mathcal{D}$ belongs to the $\alpha$-complex $\mathcal{C}_\alpha$.

If $\sigma_T$ is a tetrahedron, define $\underline{\mu}_T = \overline{\mu}_T = \varrho_T$. Otherwise,

$$\underline{\mu}_T = \min\{\varrho_{T'} \mid \sigma_{T'} \in \mathrm{up}(\sigma_T), \text{ unattached}\} \text{ and}$$
$$\overline{\mu}_T = \max\{\varrho_{T'} \mid \sigma_{T'} \in \mathrm{up}(\sigma_T)\}.$$

It is sufficient to consider only the set

$$\mathrm{up}_1(\sigma_T) = \{\sigma_{T'} \in \mathrm{up}(\sigma_T) \mid |T'| = |T| + 1\},$$

that is, all simplices incident to $\sigma_T$ whose dimension is one higher than that of $\sigma_T$, in order to derive the values $\underline{\mu}_T$ and $\overline{\mu}_T$:

$$\underline{\mu}_T = \min \left( \begin{array}{l} \{\varrho_{T'} \mid \sigma_{T'} \in \mathrm{up}_1(\sigma_T), \text{ unattached}\} \\ \cup \ \{\underline{\mu}_{T'} \mid \sigma_{T'} \in \mathrm{up}_1(\sigma_T), \text{ attached}\} \end{array} \right),$$

and

$$\overline{\mu}_T = \max\{\overline{\mu}_{T'} \mid \sigma_{T'} \in \mathrm{up}_1(\sigma_T)\}.$$

**Specifying Intervals.** The intervals of $\alpha$ values in which $\sigma_T$ is an interior, regular, or singular simplex of $\mathcal{C}_\alpha$ are shown in the table 1. Because of the general-position assumption, $\alpha$ is different from all $\varrho$ values and therefore also from all $\underline{\mu}$ and $\overline{\mu}$ values. We can thus define all intervals as open, except at endpoints 0 and $\infty$. It is necessary to distinguish simplices that bound the convex hull of $S$ from the others. The next paragraph briefly explains the entries of table 1 for the case of triangles that do not bound the convex hull of $S$. The arguments for tetrahedra, triangles on the convex hull, edges, and vertices are similar.

Consider a triangle $\sigma_T \in \mathcal{D}$, $T = \{p_i, p_j, p_k\}$, that does not bound the convex hull of $S$; we denote this by $\sigma_T \notin \partial\mathrm{conv}(S)$. Let $\sigma_{T'}$ and $\sigma_{T''}$ be the two incident tetrahedra in $\mathcal{D}$, and assume $T' = T \cup \{p_u\}$, and $T'' = T \cup \{p_v\}$. Furthermore, let $0 < \varrho_{T'} < \varrho_{T''} < \infty$; in other words, $\underline{\mu}_T = \varrho_{T'}$ and $\overline{\mu}_T = \varrho_{T''}$. Now, fix a value for $\alpha$. If $\varrho_{T''} < \alpha \leq \infty$ then the triangle $\sigma_T$ is not $\alpha$-exposed. It will, however, be part of the interior of $\mathcal{S}_\alpha$, because both incident tetrahedra are in $\mathcal{C}_\alpha$. If $\varrho_{T'} < \alpha < \varrho_{T''}$ then the triangle is $\alpha$-exposed and $\sigma_{T'}$ is in $\mathcal{C}_\alpha$ but $\sigma_{T''}$ is not. This means that $\sigma_{T'}$ is a regular triangle of $\mathcal{C}_\alpha$. For $\alpha < \varrho_{T'}$, neither $\sigma_{T'}$ nor $\sigma_{T''}$ are tetrahedra of $\mathcal{S}_\alpha$, but $\sigma_T$ can still be a singular triangle, that is, iff $\varrho_T < \alpha$ and neither $p_u$ nor $p_v$ are inside $b_T$. If one of the two points is inside $b_T$, then $\sigma_T$ is attached, and $\sigma_T$ can never be a singular triangle of $\mathcal{C}_\alpha$, no matter what $\alpha$ value is selected.

The $\alpha$-complex consists of all interior, regular, and singular simplices for a given $\alpha$ value. The interior of the $\alpha$-shape is triangulated by the interior simplices. The boundary of the interior is formed by the set of regular triangles and their edges and vertices.

Consistent with the definition in section 4, we refer to the endpoints of the intervals in table 1 as $\alpha$-thresholds. This does not include 0 and $\infty$. Since all $\underline{\mu}_T$ and $\overline{\mu}_T$ values are $\varrho$ values of other simplices, each $\alpha$-threshold is the radius of a simplex in $\mathcal{D}$. More specifically, the set of $\alpha$-thresholds is exactly the set of radii of all unattached $k$-simplices for $1 \leq k \leq 3$. Define the $\alpha$-*spectrum* as the sorted sequence of $\alpha$-thresholds. This concept will appear again in section 6.

**Computing Intervals.** Assume that each simplex $\sigma_T \in \mathcal{D}$ is marked as either "$\in \partial \text{conv}(S)$" or "$\notin \partial \text{conv}(S)$" after the construction of $\mathcal{D}$. With this, the above intervals can be computed by classifying $\sigma_T$ as attached or unattached, and by computing $\varrho_T$, $\underline{\mu}_T$, and $\overline{\mu}_T$, whenever applicable. We said that $\sigma_T$ is attached iff one of the simplices that contain $\sigma_T$ has a vertex in $b_T$, the open ball bounded by the smallest circumsphere of $\sigma_T$. This implies that $\sigma_T$ can be classified in time proportional to $|\text{up}_1(\sigma_T)|$. The time it takes to classify *all* simplices is proportional to the number of simplices in $\mathcal{D}$, because each simplex has only a constant number of faces. In other words, assuming that constant time suffices to decide whether or not a point belongs to $b_T$ (see section 5.3), a simplex can be classified in constant amortized time.

Furthermore, assume that, given $T$ with $\sigma_T \in \mathcal{D}$, $\varrho_T$ can be computed in constant time (again, see section 5.3). By processing tetrahedra before triangles before edges before vertices, we can get $\underline{\mu}_T$ and $\overline{\mu}_T$ simply as the minimum and maximum of the values $\varrho_{T'}$, $\underline{\mu}_{T'}$, and $\overline{\mu}_{T'}$, for $\sigma_{T'} \in \text{up}_1(\sigma_T)$. This also takes only constant amortized time per simplex.

## 5.3 Geometric Primitives

What are the primitive operations needed to compute $\alpha$-shapes in $\mathbb{R}^3$? Constructing Delaunay triangulations requires two geometric tests. These are a test for deciding on which side of a plane spanned by three points a fourth point lies, and one for deciding on which side of a sphere spanned by four points a fifth point lies. In order to generate the intervals of table 1, we need to compute the radius of the smallest circumsphere of a tetrahedron, triangle, or edge, and test whether a point lies inside or outside this sphere. While the two tests required for Delaunay triangulations are fairly common in geometric algorithms and computer graphics, the operations involving smallest circumspheres of triangles and edges are rather specialized. All operations share the problem of degenerate cases, which we can ignore because of the general-position assumption (see also section 6). This section gives a formula for each of the primitive operations mentioned above.

Assume that the points of $S$ are labeled as $p_1, p_2, \ldots, p_n$, and that each point $p_i$ is given by the vector $(\pi_{i,1}, \pi_{i,2}, \pi_{i,3})$ of its three coordinates. To simplify the notation in the remainder of this section, we use minors, which are determinants of submatrices of a given matrix. For convenience, define $\pi_{i,0} = 1$ for all $i$ and use the following notation for minors:

$$\mathcal{M}^{i_1,i_2,\ldots,i_k}_{j_1,j_2,\ldots,j_k} = \det \begin{pmatrix} \pi_{i_1,j_1} & \pi_{i_1,j_2} & \cdots & \pi_{i_1,j_k} \\ \pi_{i_2,j_1} & \pi_{i_2,j_2} & \cdots & \pi_{i_2,j_k} \\ \vdots & \vdots & \ddots & \vdots \\ \pi_{i_k,j_1} & \pi_{i_k,j_2} & \cdots & \pi_{i_k,j_k} \end{pmatrix}.$$

Indeed, all geometric primitives are expressed in terms of determinants, which provide a convenient and compact notation. This does not exclude the possibility of implementing the primitives using equivalent formulas; these would be, in some sense, ways of evaluating determinants that differ from the usual constructive definition.

**Plane Test.** Let $T = \{p_i, p_j, p_k\}$ and define $h_T$ as the unique plane that contains all three points of $T$. This plane can be oriented if we replace the *set* $T$ by the *sequence* $T$; for example, $T = (p_i, p_j, p_k)$. Then one side (or open half-space) of $h_T$ can be called *positive* and the other *negative*. We also refer to these as the positive and negative sides of the sequence $(p_i, p_j, p_k)$.

$$T = (p_i, p_j, p_k): \quad p_u \text{ lies on the } positive \text{ side of } h_T \iff \mathcal{M}^{i,j,k,u}_{1,2,3,0} > 0, \tag{5-1}$$

and $p_u$ lies on the *negative* side if the determinant is negative. Intuitively, $p_u$ sees the sequence of three points $p_i, p_j, p_k$ in a clockwise order iff $p_u$ lies on their positive side. Similarly, $p_u$ sees the sequence in a counterclockwise order iff $p_u$ lies on the negative side of the points. The sign of the determinant is called the *orientation* of the sequence $(p_i, p_j, p_k, p_u)$. Notice that the determinant equals zero iff the points are in degenerate position, that is, they lie on a common plane. Observe also that the orientation of a permutation of a sequence of four points is the same as the orientation of the sequence itself, provided the number of transpositions is even. Otherwise, it is the opposite. This follows trivially from the fact that the value of the determinant changes sign whenever two rows are exchanged.

**Sphere Test.** Given a set $T = \{p_i, p_j, p_k, p_u\}$, we need to decide whether another point $p_v$ lies inside or outside the sphere $\partial b_T$. We can assume that the degenerate case where $p_v \in \partial b_T$ does not occur. A possible implementation of this test is discussed in [15] using an extension of the lifting map (as mentioned in section 4) to three-dimensional spheres. Consider the paraboloid of revolution $U\colon x_4 = \sum_{\ell=1}^{3} x_\ell^2$ in $I\!\!R^4$. A sphere $\partial b$ with center $c = (\gamma_1, \gamma_2, \gamma_3)$ and radius $\rho$ is mapped to the hyperplane $\partial b_U\colon x_4 = \sum_{\ell=1}^{3}(2\gamma_\ell x_\ell - \gamma_\ell^2) + \rho^2$. This hyperplane has the property that $U \cap \partial b_U$ projected along the $x_4$-axis into the $x_1 x_2 x_3$-space yields $\partial b$. Moreover, a point $p$ lies inside (outside) $\partial b$ iff $p_U$ lies vertically below (above) $\partial b_U$. The resulting formula assumes that each point $p_i \in S$ has a fourth coordinate $\pi_{i,4} = \sum_{j=1}^{3} \pi_{i,j}^2$.

$$T = \{p_i, p_j, p_k, p_u\}: \quad p_v \text{ lies inside } \partial b_T \iff \mathcal{M}^{i,j,k,u}_{1,2,3,0} \cdot \mathcal{M}^{i,j,k,u,v}_{1,2,3,4,0} > 0 \tag{5-2}$$

The first minor, $\mathcal{M}^{i,j,k,u}_{1,2,3,0}$, is a corrective term that is necessary because the sphere does not change if the first four points are permuted. The second minor, $\mathcal{M}^{i,j,k,u,v}_{1,2,3,4,0}$, expresses the fact that the lifting map transforms a sphere test in $I\!\!R^3$ to a hyperplane test in $I\!\!R^4$.

**Radius of a Smallest Circumsphere.** Next, we consider computing the radius $\varrho_T$ of $\partial b_T$, the smallest circumsphere of $\sigma_T$, for all $k$-simplices $\sigma_T \in \mathcal{D}$, with $1 \leq k \leq 3$. The formulas for the square of $\varrho_T$ are given in (5-3) through (5-5). Note that computing $\varrho_T^2$ will be sufficient for our purposes since $\varrho_T$ can never be negative. We distinguish the cases when $k = |T| - 1$ is 1, 2, or 3.

$$T = \{p_i, p_j\}: \quad \varrho_T^2 = \frac{\left(\mathcal{M}_{1,0}^{i,j}\right)^2 + \left(\mathcal{M}_{2,0}^{i,j}\right)^2 + \left(\mathcal{M}_{3,0}^{i,j}\right)^2}{4} \tag{5-3}$$

$$T = \{p_i, p_j, p_k\}: \quad \varrho_T^2 = \frac{(\sum_{\ell=1}^{3} \left(\mathcal{M}_{\ell,0}^{i,j}\right)^2) \cdot (\sum_{\ell=1}^{3} \left(\mathcal{M}_{\ell,0}^{j,k}\right)^2) \cdot (\sum_{\ell=1}^{3} \left(\mathcal{M}_{\ell,0}^{k,i}\right)^2)}{4 \cdot (\left(\mathcal{M}_{2,3,0}^{i,j,k}\right)^2 + \left(\mathcal{M}_{1,3,0}^{i,j,k}\right)^2 + \left(\mathcal{M}_{1,2,0}^{i,j,k}\right)^2)} \tag{5-4}$$

$$T = \{p_i, p_j, p_k, p_u\}: \quad \varrho_T^2 = \frac{\left(\mathcal{M}_{2,3,4,0}^{i,j,k,u}\right)^2 + \left(\mathcal{M}_{1,3,4,0}^{i,j,k,u}\right)^2 + \left(\mathcal{M}_{1,2,4,0}^{i,j,k,u}\right)^2 + 4 \cdot \mathcal{M}_{1,2,3,0}^{i,j,k,u} \cdot \mathcal{M}_{1,2,3,4}^{i,j,k,u}}{4 \cdot \left(\mathcal{M}_{1,2,3,0}^{i,j,k,u}\right)^2} \tag{5-5}$$

In order to explain these three formulas we introduce some notation. Let $a, b, c, d$ be points in $\mathbb{R}^3$. We write $|ab|$ for the *length* of the edge $\text{conv}(\{a, b\})$, and $|abc|$ for the *area* of the triangle $\text{conv}(\{a, b, c\})$. In the case of two points, $\varrho_T$ is the same as half the distance between $p_i$ and $p_j$. Equation (5-3) follows because

$$|p_i p_j|^2 = \left(\mathcal{M}_{1,0}^{i,j}\right)^2 + \left(\mathcal{M}_{2,0}^{i,j}\right)^2 + \left(\mathcal{M}_{3,0}^{i,j}\right)^2.$$

To handle the case of a triangle, that is, $T = \{p_i, p_j, p_k\}$, we use the formulas

$$\varrho_T = \frac{|p_i p_j| \cdot |p_j p_k| \cdot |p_k p_i|}{4 \cdot |p_i p_j p_k|} \quad \text{and}$$

$$|p_i p_j p_k|^2 = \frac{1}{4} \cdot (\left(\mathcal{M}_{2,3,0}^{i,j,k}\right)^2 + \left(\mathcal{M}_{1,3,0}^{i,j,k}\right)^2 + \left(\mathcal{M}_{1,2,0}^{i,j,k}\right)^2)$$

which can be found in any good mathematical handbook. Finally, we obtain (5-5) using the extension of the lifting map mentioned above. If $\partial b$ is the sphere through points $p_i, p_j, p_k, p_u$ then $\partial b_U$ is the hyperplane through the four points $p_{i,U}, p_{j,U}, p_{k,U}, p_{u,U}$. The equation for the hyperplane can be computed directly from the coordinates of the lifted points. From this equation it is easy to compute the center and the radius of $\partial b$.

**Attached and Unattached Edges and Triangles.** We still have to consider the problem of deciding whether an edge or triangle $\sigma_T \in \mathcal{D}$ is attached or not. By definition, $\sigma_T$ is attached if there is a $\sigma_R \in \text{up}_1(\sigma_T)$ so that the point in $R - T$ belongs to $b_T$. If $\sigma_T$ is an edge, say, $T = \{p_i, p_j\}$ and $R - T = \{p_k\}$, this can be done by comparing $\varrho_T$ with the distance between $p_k$ and $\frac{p_i + p_j}{2}$. Straightforward algebraic manipulations lead to the following equation.

$$T = \{p_i, p_j\}: \quad p_k \in b_T \iff \sum_{\ell=1}^{3} \left(\mathcal{M}_{\ell,0}^{i,j}\right)^2 - \sum_{\ell=1}^{3} \left(\mathcal{M}_{\ell,0}^{i,k} + \mathcal{M}_{\ell,0}^{j,k}\right)^2 > 0 \tag{5-6}$$

Now let $\sigma_T$ be a triangle, for example, $T = \{p_i, p_j, p_k\}$ and $R - T = \{p_u\}$. To see whether or not the point $p_u$ belongs to $b_T$, we compute the center $c$ of the circumsphere $\partial b_R$ of the tetrahedron $\sigma_R$. Observe that $p_u \in b_T$ iff $c$ and $p_u$ do not lie on the same side of the plane through $p_i, p_j, p_k$. In other words, we need to test whether or not the sequences $(p_u, p_i, p_j, p_k)$ and $(c, p_i, p_j, p_k)$ have different orientation. Some rather tedious algebraic manipulations are needed to derive the following equation which expresses the derivation in terms of minors.

$$T = \{p_i, p_j, p_k\}: \ p_u \in b_T \iff \begin{cases} \mathcal{M}^{i,j,k,u}_{2,3,4,0} \cdot \mathcal{M}^{i,j,k}_{2,3,0} + \mathcal{M}^{i,j,k,u}_{1,3,4,0} \cdot \mathcal{M}^{i,j,k}_{1,3,0} \\ + \mathcal{M}^{i,j,k,u}_{1,2,4,0} \cdot \mathcal{M}^{i,j,k}_{1,2,0} - 2 \cdot \mathcal{M}^{i,j,k,u}_{1,2,3,0} \cdot \mathcal{M}^{i,j,k}_{1,2,3} > 0 \end{cases} \quad (5\text{-}7)$$

**General Position Revisited.** The general-position assumption used in this paper assures that no geometric test is ambiguous. We summarize and revise the necessary assumptions below and include pointers to the formulas for which the assumptions are relevant.

- No 4 points lie on a common plane; compare with (5-1).

- No 5 points lie on a common sphere; compare with (5-2).

- No smallest circumsphere of 2, 3, or 4 points has a radius equal to any given $\alpha$; compare with (5-3), (5-4), and (5-5).

- No point lies on the smallest circumsphere of 2 or 3 other points; compare with (5-6) and (5-7).

These assumptions are indeed very restrictive and rarely true for real-life data. We will deal with this apparent shortcoming in section 6.2.

# 6 Implementation

Our current implementation of a software tool for $\alpha$-shapes in $\mathbb{R}^3$ consists of the following three parts.

1. A program that constructs Delaunay triangulations using flips (see section 5.1).

2. A program that computes the $\alpha$-intervals for all simplices in a Delaunay triangulation, and then sorts the endpoints of these intervals (see section 5.2).

3. An $\alpha$-shape visualizer that enables the user to manually select different $\alpha$ values and render the corresponding shape on a graphics workstation (see figures 1 through 5).

Parts 1 and 2 are preprocessing steps that take time $O(n^2)$ and $O(m \log m)$, where $n$ is the number of points and $m$ is the number of simplices of $\mathcal{D}$. The current code for part 3 takes time $O(m)$ to render a particular $\alpha$-shape. Improvements based on fast data structures for intervals are forthcoming.

One important aspect of the implementation is its *robustness*. By this we mean that the program will either produce the correct output for the original data set $S$, or, in case of degeneracies, it guarantees to give the correct result for a set $S(\varepsilon)$ arbitrarily close to the original input. This is achieved by a symbolic perturbation scheme briefly described in section 6.2. The method avoids possible conflicts between the topological and geometric structure of the data by using exact arithmetic and by perturbing the original data set such that all degeneracies disappear. The program can thus be considered to be purely combinatorial and logical so that correctness in the strict sense is possible in principle. Note that the perturbation is infinitesimal and symbolic, that is, $S$ and $S(\varepsilon)$ can be viewed as arbitrarily close together, and the computational overhead is independent of $\varepsilon$.

## 6.1 Data Structures

There are two main data structures needed for storing the family of $\alpha$-shapes of a given data set. One represents the connectivity and order among the simplices of the three-dimensional Delaunay triangulations. The other is used for the intervals assigned to the simplices of $\mathcal{D}$. A triangle-based data structure is used to store $\mathcal{D}$. This is briefly described in the paragraphs below. An interval tree can be used to store the collection of intervals (see for example [35]). The current version of our program, however, stores the $\alpha$-spectrum using only a linear array. Recall that the $\alpha$-spectrum is the sorted sequence of $\alpha$-thresholds, and for practical reasons, the implementation adds 0 and $\infty$ to this list. Universal hashing (see for example [5]) provides the link between the triangle structure and the array.

The data structure used to store the three-dimensional triangulation is a specialized version of the edge-facet structure introduced by Dobkin and Laszlo [7]. Related data structures are the quad-edge structure due to Guibas and Stolfi [22], which can be used to model two-dimensional manifolds, and the cell-tuple structure by Brisson [2], which works in arbitrary dimensions. The edge-facet structure is designed for general cell complexes in three dimensions. By reducing the scope to simplicial complexes, it is possible to improve the compactness and the speed of the structure. We refer to the result as the *triangle-edge structure*.

The atomic unit of the triangle-edge structure is the so-called *triangle-edge pair* $a = \langle \sigma, i \rangle$, with $0 \leq i \leq 5$. It identifies six versions of the triangle $\sigma$, one for each of its six directed edges. Each triangle defines two *edge rings*. One edge ring traverses the edges of $\sigma$ in a counterclockwise order, the other traverses them in a clockwise order. Similarly, each edge defines two *triangle rings* traversing the incident triangles in the two opposite orders. Each triangle-edge $a$ belongs to exactly one edge ring and exactly one triangle ring.

The internal representation of the structure takes advantage of the fact that all edge rings have length three. It is thus possible to avoid actual pointers for the edge rings by merging the six triangle-edge pairs of two opposite edge rings into one record. Such a record allocates 30 (36) bytes per triangle, assuming that 2-byte (4-byte) integers are used as indices to the vertices, and 4-byte integers for triangle-edge pairs. Further details are omitted.

## 6.2 Simulated Perturbation

For implementation purposes it is no longer appropriate to assume that the input points are in general position. This assumption would be too restrictive. On the other hand, in the context of three-dimensional $\alpha$-shapes, it would be rather tedious to deal with the large number of special cases in ad hoc manner. For this reason, we apply a general technique, known as Simulation of Simplicity, or SoS [15], which acts as a black box between the implementation of a geometric algorithm and the input data. It allows a systematic treatment of all special cases on the level of geometric primitive operations. The SoS library consists of a set of carefully implemented primitives. It provides the programmer with the illusion of simple data while the actual input is in arbitrary and thus possibly degenerate position.[8] This section can only sketch the basic idea of SoS (refer to [15] for details).

The idea of SoS is to perturb the given objects ever so slightly, in a manner that *all* degeneracies disappear. The perturbation should be small enough so that the nondegenerate position of objects relative to each other remains unchanged. Since it is usually rather difficult and costly to actually come up with such a perturbation, SoS performs it symbolically by substituting polynomials in $\varepsilon$ for the parameters specifying each object. In the context of this paper, the coordinate $\pi_{i,j}$ of the input point $p_i \in S$ is replaced by its *perturbed* version

$$\pi_{i,j}(\varepsilon) = \pi_{i,j} + \varepsilon(i,j),$$

where $\varepsilon(i,j) = \varepsilon^{\delta^{4i-j}}$ and $\delta$ is a sufficiently large constant. The choice of $\varepsilon(i,j)$ implies that general position of the perturbed input set $S(\varepsilon)$ is guaranteed if $\varepsilon$ is positive and arbitrarily small. Thus, $\varepsilon$ can be treated as an indeterminant in the symbolic evaluation of the primitive operations which are now based on $S(\varepsilon)$.

The expressions that correspond to primitive operations are polynomials in $\varepsilon$. The coefficients of these polynomials are expressions in the coordinates of the original data. For example, in order to determine whether the sequence $(p_i, p_j, p_k, p_u)$, with $i < j < k < u$, has positive orientation, the sign of the polynomial

$$\mathcal{M}_{1,2,3,0}^{i,j,k,u} + \mathcal{M}_{1,2,0}^{j,k,u} \cdot \varepsilon(i,3) - \mathcal{M}_{1,3,0}^{j,k,u} \cdot \varepsilon(i,2) + \mathcal{M}_{2,3,1}^{j,k,u} \cdot \varepsilon(i,1) + \cdots$$

$$\cdots + \mathcal{M}_{2,0}^{j,u} \cdot \varepsilon(k,3)\varepsilon(i,1) + \mathcal{M}_{1,0}^{i,u} \cdot \varepsilon(k,3)\varepsilon(j,2) + \varepsilon(k,3)\varepsilon(j,2)\varepsilon(i,1) + \cdots$$

must be evaluated. Assume that the terms in the polynomial are sorted in order of increasing exponents of $\varepsilon$. We say the evaluation has *depth* $t$ if the coefficient of the $(t+1)$-st term is the first one that does not vanish. Because $\varepsilon$ is assumed to be arbitrarily small, this term decides the sign of the entire polynomial. Notice, that the coefficient at depth 0 is the same as the expression of the primitive for the original data. This implies that $S(\varepsilon)$ is consistent with $S$ as far as nondegenerate configurations are concerned. Observe also that the polynomial never evaluates to zero since there is always a term with non-zero coefficient. In other words, $S(\varepsilon)$ is simple.

---

[8]The terms "simple," "general position," and "nondegenerate position" are used as synonyms. Notice how the "general case" of the algorithm designer is usually the "simple case" for the implementing programmer.

The code that implements the polynomial of a given primitive operation can be generated automatically. The overhead for the symbolic perturbation as such is thus neglectable. However, the perturbation can only be simulated when exact arithmetic is used to compute the coefficients of the polynomials. More precisely, we need to be able to tell when a coefficient vanishes. Of course, exact arithmetic entails a somewhat higher cost for the low-level computations, but we believe that this is an adequate price for a compact and robust implementation of a rather involved geometric algorithm.

Note that SoS does not allow the user to selectively remove certain types of degeneracies, while others remain. Rather, it *fully* implements what theoreticians call "general position." Observe that algorithms employing SoS produce results for the perturbed data, and not for the original one. Some postprocessing can be used to remove or repair parts of the output that collapse to a degenerate state because of the degeneracy of the input data. This is mentioned in [15]. As an example, consider the case of Delaunay triangulations and coplanar points on the boundary of the convex hull. The perturbation will move some points towards inside and some towards outside. In the Delaunay triangulation, the points that moved from the convex hull boundary inside will be covered by flat tetrahedra. These tetrahedra have infinitesimal thickness in the perturbed setting, and are completely flat in the original setting. It is easy to identify and remove these tetrahedra in a postprocessing step.

### 6.3 Performance

Table 2 shows the performance of the incremental-flip algorithm (see section 5.1) for a number of test runs. We count the number of flips, and the number of necessary evaluations of 5-by-5 and 4-by-4 determinants. The "max" and "mean" depth columns count the "depths" of the corresponding evaluations (see section 6.2). These columns give a measure of *how* degenerated a data set is. The experiments were run on Silicon Graphics workstations with 50 MHz MIPS R4000 CPUs and 48 Mb or more of main memory. Sample frames for the data sets `molecule`, `tori`, `universe`, `bust`, and `phone_1` can be found in figures 1 through 5.

The running time for the presented examples is way better than predicted by the quadratic worst case. This is due to the fact that Delaunay triangulations usually do not reach their worst-case upper bounds of $\Theta(n^2)$ faces (see section 4). Table 2 shows running times that seem to be roughly proportional to $n(\log n)^2$. An exception to this rule are volumetric data with points on a regular grid or parallel slices (as in the data sets `grid1`, `ku2`, `rat_T2a`, and `hsr`). The incremental-flip algorithm is certainly less than ideal for such distributions and should possibly be replaced by a randomized version (see [16]).[9]

Observe the positive correlation between the number of determinant evaluations or long-integer operations and CPU seconds. Indeed, run-time profiles of the `C` code suggest that the majority of CPU cycles is used for the long-integer arithmetic computing determinants. We estimate that approximately 75% of the time is spent on long-integer arithmetic. The multiplication routine is responsible for more than half of it.

Table 3 shows the performance of the program that generates the $\alpha$-intervals for all simplices in the Delaunay triangulation (see section 5.2). The number of simplices in a triangulation is denoted by $|\mathcal{D}| = \sum_{k=1}^{3} |F_k|$. A large portion of the memory requirement is due to the fact that, in order to achieve robustness, we need to compute the $\alpha$-thresholds with long-integer arithmetic. The integers involved can get fairly long because the corresponding expressions are rather involved (see expressions (5-3) through (5-5) in section 5.3). However, as soon as their correct order is determined, the exact values of the $\alpha$-thresholds, which are long-integer rationals, are no longer needed. For rendering purposes, standard floating-point accuracy is certainly sufficient. Duplicate $\alpha$-thresholds, including the ones caused by attached simplices, are eliminated in the sorting phase. The size of the $\alpha$-spectrum, that is, the final number of $\alpha$-thresholds, plus the values 0 and $\infty$, is usually considerably smaller than the number of simplices in the triangulation.

---

[9]As a matter of fact, the second author implemented a variant of the randomized algorithm, achieving performance numbers roughly twice as fast than the ones in table 2 (see [34]).

| $S$ | $|S|$ | $|F_{0,\infty}|$ | $\mathcal{D}$ | | flips | SoS primitives | | SoS depth | | long-integer arithmetic | | memory | CPU |
|---|---|---|---|---|---|---|---|---|---|---|---|---|---|
| | | | $|F_2|$ | $|F_3|$ | | $\mathcal{M}^{i,j,k,u,v}_{1,2,3,4,0}$ | $\mathcal{M}^{i,j,k,u}_{1,2,3,0}$ | max | mean | $+/-$ | $*$ | Mb | seconds |
| molecule | 318 | 34 | 4000 | 1984 | 7117 | 43745 | 21823 | 3 | 0.000061 | 1755229 | 1428355 | 0.383 | 10.38 |
| ma1 | 450 | 120 | 5346 | 2614 | 6299 | 39691 | 22221 | 12 | 0.031787 | 1600911 | 1330478 | 0.516 | 8.97 |
| tori | 800 | 193 | 12197 | 6003 | 19483 | 120227 | 61114 | 0 | 0.000000 | 4855622 | 3937343 | 1.789 | 30.94 |
| grid1 | 1000 | 272 | 12192 | 5961 | 42744 | 257787 | 122940 | 13 | 0.892645 | 10494784 | 10387794 | 1.125 | 61.33 |
| spiral | 1248 | 191 | 19301 | 9556 | 49329 | 301479 | 150786 | 13 | 0.098697 | 11927947 | 10058019 | 2.773 | 74.93 |
| molecule2 | 1366 | 70 | 17974 | 8953 | 44908 | 274912 | 129670 | 0 | 0.000000 | 10945857 | 8903536 | 1.523 | 68.80 |
| universe | 1717 | 77 | 21321 | 10623 | 49788 | 305928 | 147307 | 0 | 0.000000 | 12223778 | 9942300 | 1.852 | 68.71 |
| pdb9pap | 1908 | 54 | 25318 | 12633 | 64607 | 395586 | 184959 | 0 | 0.000000 | 15755170 | 12825709 | 2.195 | 92.58 |
| ma2 | 2206 | 871 | 46427 | 22779 | 139984 | 867839 | 447606 | 4 | 0.017762 | 34764277 | 28516841 | 5.328 | 206.77 |
| bust | 2630 | 248 | 35196 | 17475 | 101267 | 618854 | 291070 | 8 | 0.090630 | 24605560 | 20488777 | 2.906 | 160.13 |
| fract2 | 2704 | 141 | 35203 | 17532 | 181432 | 1098618 | 524976 | 8 | 0.127745 | 42864764 | 36889289 | 2.953 | 272.55 |
| br | 3762 | 89 | 51559 | 25736 | 142562 | 871651 | 401509 | 3 | 0.000002 | 34550899 | 28142998 | 4.406 | 201.92 |
| liss_5_8 | 4200 | 534 | 58486 | 28977 | 142717 | 872057 | 425637 | 4 | 0.000349 | 34933465 | 28389761 | 4.961 | 215.78 |
| teapot | 4668 | 1247 | 55529 | 27142 | 178761 | 1094305 | 552003 | 9 | 0.056926 | 43311204 | 36476581 | 4.914 | 277.65 |
| ku2 | 5292 | 462 | 74836 | 37188 | 1371107 | 8218910 | 3779546 | 13 | 0.345176 | 322279363 | 291020300 | 5.977 | 1855.50 |
| phone_1 | 6070 | 105 | 80131 | 40014 | 251224 | 1534163 | 705118 | 13 | 0.013700 | 59982685 | 49669888 | 6.969 | 407.13 |
| rat_T2a | 9600 | 263 | 125933 | 62836 | 1170671 | 7089183 | 3480820 | 11 | 1.717422 | 386072493 | 352580891 | 10.984 | 2096.14 |
| ra | 10000 | 566 | 130172 | 64804 | 412360 | 2514540 | 1171071 | 2 | 0.000008 | 99914656 | 81317780 | 11.391 | 627.12 |
| hsr | 10088 | 131 | 122077 | 60974 | 942250 | 5694210 | 2586581 | 13 | 0.758750 | 241479463 | 226164544 | 10.664 | 1369.19 |
| smoke015 | 15000 | 672 | 192898 | 96114 | 636703 | 3881595 | 1801968 | 0 | 0.000000 | 154212466 | 125483265 | 16.969 | 1022.71 |

Table 2: Performance of the incremental-flip algorithm.

| $S$ | $|S|$ | $|\mathcal{D}|$ | |spectrum| | SoS primitives | | | SoS depth | | long-integer arithmetic | | memory | CPU |
|---|---|---|---|---|---|---|---|---|---|---|---|---|
| | | | | Expr (5-5) | Expr (5-4) | Expr (5-3) | max | mean | $+/-$ | $*$ | Mb | seconds |
| molecule | 318 | 8317 | 4067 | 1984 | 1469 | 632 | 0 | 0.000000 | 956974 | 732611 | 0.930 | 7.80 |
| ma1 | 450 | 11141 | 5164 | 2614 | 1754 | 1071 | 37 | 0.016009 | 1232528 | 968669 | 1.172 | 9.37 |
| tori | 800 | 25193 | 7811 | 6003 | 3338 | 1848 | 0 | 0.000000 | 2808806 | 2153544 | 2.680 | 27.23 |
| grid1 | 1000 | 25383 | 5 | 5961 | 4854 | 2700 | 39 | 0.490012 | 2141341 | 1911438 | 2.562 | 16.12 |
| spiral | 1248 | 39849 | 7939 | 9556 | 4710 | 2649 | 37 | 0.009140 | 4289522 | 3357255 | 4.195 | 40.87 |
| molecule2 | 1366 | 37313 | 20475 | 8953 | 8050 | 3473 | 0 | 0.000000 | 4573791 | 3592201 | 4.141 | 44.90 |
| universe | 1717 | 44358 | 21443 | 10623 | 7155 | 3690 | 0 | 0.000000 | 5166676 | 4077980 | 4.539 | 41.19 |
| pdb9pap | 1908 | 52543 | 31146 | 12633 | 12515 | 5998 | 0 | 0.000000 | 5298887 | 6665298 | 5.672 | 58.07 |
| ma2 | 2206 | 95056 | 32434 | 22779 | 6611 | 3046 | 0 | 0.000000 | 10081520 | 7799049 | 19.078 | 82.66 |
| bust | 2630 | 73021 | 34041 | 17475 | 11799 | 7541 | 2 | 0.006083 | 8408574 | 6776499 | 8.023 | 90.05 |
| fract2 | 2704 | 73109 | 31336 | 17532 | 9399 | 7465 | 24 | 0.044352 | 8010544 | 6575623 | 7.812 | 83.62 |
| br | 3762 | 106879 | 61552 | 25736 | 24711 | 11121 | 0 | 0.000000 | 13511423 | 10827946 | 11.383 | 116.29 |
| liss_5_8 | 4200 | 121171 | 56730 | 28977 | 17121 | 11943 | 0 | 0.000000 | 14056712 | 11096304 | 12.781 | 139.64 |
| teapot | 4668 | 115725 | 15434 | 27142 | 16626 | 10366 | 0 | 0.000000 | 12494168 | 10167622 | 12.164 | 128.27 |
| ku2 | 5292 | 154963 | 60715 | 37188 | 33551 | 19378 | 39 | 0.093015 | 18000791 | 15529705 | 15.977 | 153.95 |
| phone_1 | 6070 | 166331 | 80590 | 40014 | 25444 | 19773 | 24 | 0.000936 | 19477178 | 15850522 | 18.234 | 227.27 |
| rat_T2a | 9600 | 261465 | 36543 | 62836 | 56180 | 28668 | 24 | 0.016854 | 30276256 | 24180721 | 25.812 | 237.65 |
| ra | 10000 | 270343 | 132353 | 64804 | 42307 | 25245 | 2 | 0.000042 | 32007436 | 25724425 | 28.445 | 331.43 |
| hsr | 10088 | 254241 | 5763 | 60974 | 49469 | 24431 | 21 | 0.003192 | 23493214 | 18959384 | 24.836 | 186.20 |
| smoke015 | 15000 | 400795 | 200163 | 96114 | 65121 | 38926 | 0 | 0.000000 | 47929411 | 38659261 | 43.898 | 594.94 |

Table 3: Performance of the $\alpha$-interval generator.

# 7 Applications and Further Illustrations

It is important to point out that $\alpha$-shapes are a fairly generic tool that can be used in many applications that have to do with shape, including automatic mesh generation and geometric modeling (see figure 3). Indeed, they can be used as a concrete expression of shape, which is often all that is needed. Similarly, three-dimensional $\alpha$-shapes can be used to identify clusters in trivariate data. Beyond these generic applications, there are others that rely on particular properties of $\alpha$-shapes. For these applications, it would be difficult to replace $\alpha$-shapes by any other reasonable notion of shape. Two such applications are briefly addressed.

**Molecular Structures.** Molecules are usually modeled as conglomerates of atoms with fixed relative positions. Each atom is represented by a ball around a center point, and the radius of the ball depends on what the model is supposed to express. For example, in the so-called space-filling diagram (see section 3.1) the balls encompass the idealized locations of the electrons so that balls of nearby atoms typically overlap. This diagram, defined as a union of balls, is in a strict geometric sense dual to the $\alpha$-shape of the center points, assuming each ball has radius equal to $\alpha$. The $\alpha$-shape can thus be used to compute structural properties of the space-filling diagram, such as its connectivity. Alternatively, the $\alpha$-shape itself, for this value of $\alpha$, can be used to model and manipulate the molecule. When different atoms are represented by balls of *different* sizes then weighted $\alpha$-shapes need to be used (see section 3.4).

Molecules with intersting $\alpha$-shapes arise in the study of proteins and how they fold (see figure 4). The geometric locations of the atoms of about 500 proteins are known today. However, there are many more gene sequences that can be determined. One of the goals of theoretical molecular biology is to obtain three-dimensional positional information from the knowledge of these sequences. This is the problem of protein folding [20, 36]. Since the $\alpha$-shape is computationally inexpensive and because it closely reflects the physical reality of molecules, it is hoped that $\alpha$-shapes prove to be a useful tool in future protein folding research.

**Distribution of a Point Set.** An interesting though ill-defined geometric problem arises in the study of the distribution of galaxies in our universe. As observed in studies like in [4, 19], the galaxies are distributed in an unexpected and rather nonuniform manner (see figure 5). Astronomers have measured the location of about 170000 galaxies, each one represented by a point in three-dimensional space. It appears that a large number of galaxies are located on or close to sheet-like and to filament-shaped structures. In other words, large subsets of the points are distributed in a predominantly two- or one-dimensional manner.

How can this intuitive notion of the dimension of a point distribution be captured? A possible answer can be given by considering the entire family of $\alpha$-shapes. Let $A(\alpha)$ be the surface area of the $\alpha$-shape and let $V(\alpha)$ be its volume. To measure the degree to which the points are two-dimensionally distributed, it might be interesting to investigate the relative variation of $A$ and $V$ over the range of $\alpha$ values between 0 and $\infty$. More generally, it would be interesting to study the relationship between $\alpha$-shapes ant the notion of fractals and fractal dimension. Through the availability of signatures of measures (see section 8), the $\alpha$-shape gives comparably efficient access to metric information over a range of detail monitored by $\alpha$. Resulting quantifications can be useful in the comparative study of the actually observed galaxy distribution and simulated data (see [10] for first steps in this direction).

Figure 3: Phone. [$n = 6070$, $|F_2| = 80131$]

The points are obtained from a public domain data set for modeling and rendering programs. All connectivity information (edges and triangles) is removed, and in order to generate more points, the centroids of all triangles are added. As described in section 3, each $\alpha$-shape is triangulated by the tetrahedra of the corresponding $\alpha$-complex. This might be useful in the automatic generation of meshes for objects with nonconvex surfaces.

Figure 4: Molecule. [$n = 318$, $|F_2| = 4000$]
The data represents a time-averaged molecular dynamics structure of gramicidin A, a peptide that forms a channel for ion and water movement across lipid membranes. The major structural motif is a right-handed beta-bonded helix. The tunnel of the macro-structure can be detected using relatively large $\alpha$ values. Smaller values of $\alpha$ result in $\alpha$-shapes with larger numbers of isolated triangles and edges. These $\alpha$-shapes reveal the helix of the micro-structure.

Figure 5: Universe. $[n = 1717, |F_2| = 21321]$

This data represents a simulation of the positions of galaxies within our universe. The theory is that galaxies first clustered into sheet-like structures, then progressed to filament-shaped structures at the intersection of multiple sheets. As filaments began to intersect, global clusters appeared. It is interesting to investigate the macro- and micro-structure of the galaxies, including the detection of large voids and local or global clusters. The full spectrum of $\alpha$-shapes promises to be useful in this study.

# 8  Summary and Open Problems

The main contribution of this paper is the introduction of a sound framework that formally captures the rather intuitive notion of the "shape" of a point set in space. This is the concept of three-dimensional $\alpha$-shapes. A prototype version of a robust $\alpha$-shape tool has been implemented. The authors of this paper hope that this tool will find many users within the engineering and the scientific computing and visualization communities. However, there is still a lot of work to be done. For example, the extensions mentioned in section 3.4 are worthwhile implementing, and this is planned in the near future. The extensions mentioned below are either less specific or theoretically not well understood.

**Improving the Running Time.** A large fraction of the time used to construct $\alpha$-shapes is needed for computing the Delaunay triangulation of the points. The algorithm used in our implementation takes time $O(n^2)$ in the worst case, independent of the number of simplices of $\mathcal{D}$. However, it rarely exhibits worst-case behavior. Still, it would be useful to have an algorithm whose running time is roughly proportional to the size of $\mathcal{D}$. Is it possible to construct $\mathcal{D}$ in time $O(n \log n + m)$, where $m = |F_1 \cup F_2 \cup F_3|$? A first step towards such an algorithm is the output-sensitive convex hull algorithm of Seidel [37]. If combined with the methods of Matoušek and Schwarzkopf [30] it runs in time $O(n^{4/3+\varepsilon} + m \log n)$ for $n$ points in $\mathbb{R}^3$. By randomization, one can also achieve an expected running time roughly proportional to the expected number of simplices, if an underlying distribution is assumed. This is explained in [16].

On a different level, the running time of our program can be improved by speeding up the geometric primitives which all reduce to integer computations (see section 6). According to our experimental studies, about 75% of the time is spent doing integer arithmetic. This implies that appropriate hardware support might have a significant impact on the running time.

**Maintaining Alpha Shapes.** In some applications it is necessary to construct $\alpha$-shapes across a number of different point sets, and often these point sets are very similar to each other. For example, a point set might undergo local changes within an iterative process, and the $\alpha$-shape or some feature of it is to be constructed at each step of the iteration. A local change might be the insertion of a new point, the deletion of an old point, the dislocation of one point, the dislocation of a subset of the points, etc. The development of efficient update algorithms that reuse available structure as much as possible can lead to dramatic improvements of the overall running time.

**Features and Signatures.** Individual $\alpha$-shapes are interesting geometric objects, and it would be useful to have efficient algorithms that can analyze its geometric and topological properties or features. For example, computing the volume is fairly straightforward because the $\alpha$-complex provides a triangulation of the $\alpha$-shape, and the volume of the $\alpha$-shape is simply the sum of the volumes of the tetrahedra. More challenging is the computation of the connectivity of the $\alpha$-shape as expressed by its first three homology groups (see for example Giblin [21]).

As suggested by the second application discussed in section 7, the history of a feature, over all values of $\alpha$ from 0 through $\infty$, is of interest. Consider some specific feature, say, the number of connected components of the $\alpha$-shape. The corresponding *signature* is a function $c: [0, +\infty] \to \mathbb{R}$, so that $c(\alpha)$ is the number of components of $\mathcal{S}_\alpha$. This function reflects the evolution of the number of components as $\alpha$ changes continuously from 0 to $+\infty$. Given the $\alpha$-spectrum, it is fairly easy to compute $c$. Start at $\alpha = 0$ and maintain a union-find data structure (see for example [5]) storing the components as threshold values are processed in increasing order. A more challenging task is the computation of the signatures for higher-order homology groups. Such signatures might be handy in the selection of an appropriate $\alpha$ value which typically depends on the application and the user's momentary focus of attention.

# Available Software

The implementation mentioned in section 6 is called "Alvis — A 3D Alpha Shape Visualizer." This tool allows the user to interactively select $\alpha$-values and display the corresponding shape. A small collection of signatures aids the selection process. We see the program as an extension to the paper, as an "animated figure," so to speak. One of its purposes is to effectively convey the concept of $\alpha$-shapes to the engeneering and scientific computing community. It is available via anonymous ftp from `ftp.ncsa.uiuc.edu (141.142.20.50)`. The latest release of the code, including the code for three-dimensional Delaunay triangulations, can be found in the directory `SGI/Alpha-shape`. The code requires a Silicon Graphics workstation, running Irix 4.0, or later; 32 Mb main memory, or more, are advisable. This code is a new kind of shareware: you share with us your experience in applying Alvis to engineering and science problems, and we do our best to develop the software so it can meet your needs. For questions contact `<alpha@ncsa.uiuc.edu>`.

# Acknowledgements

We thank Ping Fu, for her interest in the $\alpha$-shape tool and her active encouragement and collaboration. We are also indebted to Marc Dyksterhouse, Eric Jakobson, Patrick Moran, Michael Norman, and Shankar Subramanian from NCSA for their valuable comments and interest in our work, and the staff at the NCSA Renaissance Experimental Lab for granting us access to their graphics workstations. The second author owes Donald Hearn from the CS Department at UIUC for triggering his interest in implementing an $\alpha$-shape tool with a course on scientific visualization.

"Don't look like a convex hull ... get yourself in $\alpha$-shape!"
— David Knapp.

# References

[1] J D Boissonnat. Geometric structures for three-dimensional shape representation. *ACM Transactions on Graphics*, 3(4):266–286, 1984.

[2] E Brisson. Representing geometric structures in *d* dimensions: Topology and order. *Discrete and Computational Geometry*, 9(4):387–426, 1993.

[3] A Brønsted. *An Introduction to Convex Polytopes*. Graduate Texts in Mathematics. Springer-Verlag, New York, 1983.

[4] R Y Cen, A Jameson, F Liu, and J P Ostriker. The universe in a box: Thermal effects in a standard cold dark matter scenario. *Astrophysical Journal (Letters)*, 362:L41, 1990.

[5] T H Cormen, C E Leiserson, and R L Rivest. *Introduction to Algorithms*. MIT Press, Cambridge, Massachusetts, 1990.

[6] B Delaunay. Sur la sphère vide. *Izvestia Akademii Nauk SSSR, Otdelenie Matematicheskii i Estestvennyka Nauk*, 7:793–800, 1934.

[7] D P Dobkin and M J Laszlo. Primitives for the manipulation of three-dimensional subdivisions. *Algorithmica*, 4(1):3–32, 1989.

[8] R A Drebin, L Carpenter, and P Hanrahan. Volume rendering. *Computer Graphics*, 22(4):65–74, 1988. ACM Siggraph '88. Conference Proceedings.

[9] R A Dwyer. Higher-dimensional Voronoi diagrams in linear expected time. *Discrete and Computational Geometry*, 6(4):343–367, 1991.

[10] M D Dyksterhouse. An alpha-shape view of our universe. Master's thesis, Department of Computer Science, University of Illinois at Urbana-Champaign, Urbana, Illinois, 1992.

[11] H Edelsbrunner. *Algorithms in Combinatorial Geometry*. Springer-Verlag, Berlin, 1987.

[12] H Edelsbrunner. Geometric algorithms. In P Gruber and J Wills, editors, *Handbook of Convex Geometry*, pages 699–735. North-Holland, Amsterdam, 1992.

[13] H Edelsbrunner. Weighted alpha shapes. Technical Report UIUCDCS-R-92-1760, Department of Computer Science, University of Illinois at Urbana-Champaign, Urbana, Illinois, 1992.

[14] H Edelsbrunner, D G Kirkpatrick, and R Seidel. On the shape of a set of points in the plane. *IEEE Transactions on Information Theory*, IT-29(4):551–559, 1983.

[15] H Edelsbrunner and E P Mücke. Simulation of Simplicity: A technique to cope with degenerate cases in geometric algorithms. *ACM Transactions on Graphics*, 9(1):66–104, 1990.

[16] H Edelsbrunner and N R Shah. Incremental topological flipping works for regular triangulations. In *Proceedings of the Eighth Annual Symposium on Computational Geometry*, pages 43–52, 1992.

[17] J Fairfield. Contoured shape generation: Forms that people see in dot patterns. In *Proceedings of the IEEE Conference on Cybernetics and Society*, pages 60–64, 1979.

[18] J Fairfield. Segmenting dot patterns by Voronoi diagram concavity. *IEEE Transactions on Pattern Analysis and Machine Intelligence*, 5(1):104–110, 1983.

[19] M J Geller and J P Huchra. Mapping the universe. *Science*, 246:897–903, 1989.

[20] C Ghélis and J Yon. *Protein Folding*. Academic Press, 1982.

[21] P J Giblin. *Graphs, Surfaces, and Homology*. Second edition. Chapman and Hall, London, 1981.

[22] L J Guibas and J Stolfi. Primitives for manipulation of general subdivisions and the computation of Voronoi diagrams. *ACM Transactions on Graphics*, 4(2):74–123, 1985.

[23] R A Jarvis. Computing the shape hull of points in the plane. In *Proceedings of the IEEE Computing Society Conference on Pattern Recognition and Image Processing*, pages 231–241, 1977.

[24] B Joe. Three-dimensional triangulations from local transformations. *SIAM Journal on Scientific and Statistical Computing*, 10(4):718–741, 1989.

[25] B Joe. Construction of three-dimensional Delaunay triangulations using local transformations. *Computer Aided Geometric Design*, 8(2):123–142, 1991.

[26] D G Kirkpatrick and J D Radke. A framework for computational morphology. In G T Toussaint, editor, *Computational Geometry*, pages 234–244. Elsevier North Holland, New York, 1985.

[27] C L Lawson. Software for $C^1$ surface interpolation. In J R Rice, editor, *Mathematical Software III*, pages 161–194. Academic Press, New York, 1977.

[28] C Lee. Regular triangulations of convex polytopes. In P Gritzmann and B Sturmfels, editors, *Applied Geometry and Discrete Mathematics: The Victor Klee Festschrift*, pages 443–456. American Mathematical Society, Providence, RI, 1991.

[29] W E Lorensen and H E Cline. Marching cubes: A high resolution 3D surface construction algorithm. *Computer Graphics*, 21(4):163–169, 1987. ACM Siggraph '87. Conference Proceedings.

[30] J Matoušek and O Schwarzkopf. Linear optimization queries. In *Proceedings of the Eighth Annual Symposium on Computational Geometry*, pages 16–25, 1992.

[31] D W Matula and R R Sokal. Properties of Gabriel graphs relevant to geographic variation research and the clustering of points in the plane. *Geographical Analysis*, 12(3):205–222, 1980.

[32] P McMullen and G C Shepard. *Convex Polytopes and the Upper Bound Conjecture*. Cambridge Unversity Press, London, 1971.

[33] W W Moss. Some new analytic and graphic approaches to numerical taxonomy, with an example from the Dermanyssidae (Acari). *Systematic Zoology*, 16:177–207, 1967.

[34] E P Mücke. *Shapes and Implementations in Three-Dimensional Geometry*. PhD thesis, Department of Computer Science, University of Illinois at Urbana-Champaign, Ubana, Illinois, 1993. Also as Technical Report UIUCDCS-R-93-1836.

[35] F P Preparata and M I Shamos. *Computational Geometry — An Introduction*. Springer-Verlag, New York, 1985.

[36] F M Richards. The protein folding problem. *Scientific American*, 264(1):54–63, 1991.

[37] R Seidel. Constructing higher-dimensional convex hulls in logarithmic cost per face. In *Proceedings of the 18th Annual ACM Symposium on Theory of Computing*, pages 484–507, 1986.

[38] R Seidel. Exact upper bounds for the number of faces in $d$-dimensional Voronoi diagrams. In P Gritzmann and B Sturmfels, editors, *Applied Geometry and Discrete Mathematics: The Victor Klee Festschrift*, pages 517–530. American Mathematical Society, Providence, RI, 1991.

[39] G T Toussaint. The relative neighborhood graph of a finite planar set. *Pattern Recognition*, 12(4):261–268, 1980.

[40] R C Veltkamp. *Closed Object Boundaries from Scattered Points*. PhD thesis, Erasmus Universiteit, Rotterdam, 1992.

[41] G Voronoi. Nouvelles applications des paramètres continus à la théorie des formes quadratiques. Premier Mémoire: Sur quelques propriétés des formes quadratiques positives parfaites. *Journal für die Reine und Angewandte Mathematik*, 133:97–178, 1907.

[42] G Voronoi. Nouvelles applications des paramètres continus à la théorie des formes quadratiques. Deuxième Mémoire: Recherches sur les paralléllöedres primitifs. *Journal für die Reine und Angewandte Mathematik*, 134:198–287, 1908.